\documentclass[11pt]{amsart}
\usepackage{amssymb,amsmath}
\usepackage{amsthm,enumerate}
\usepackage[applemac]{inputenc}

\newtheorem{thm}[equation]{Theorem}
\newtheorem{lem}[equation]{Lemma}
\newtheorem{prop}[equation]{Proposition}

\newtheorem{rem}[equation]{Remark}
\numberwithin{equation}{section}

\newcommand{\R}{\mathbb{R}}
\newcommand{\N}{\mathbb{N}}
\renewcommand{\H}{\mathrm{H}}
\newcommand{\T}{\mathbf{T}}

\begin{document}
\title[Transport-Entropy inequalities on the line]
{Transport-Entropy inequalities on the line}
\author{Nathael Gozlan}

\date{\today}

\address{Universit\'e Paris Est Marne la Vall\'ee - Laboratoire d'Analyse et de Math\'e\-matiques Appliqu\'ees (UMR CNRS 8050), 5 bd Descartes, 77454 Marne la Vall\'ee Cedex 2, France}
\email{nathael.gozlan@univ-mlv.fr}


\maketitle

\begin{abstract}
We give a necessary and sufficient condition for transport-entropy inequalities in dimension one. As an application, we construct a new example of a probability distribution verifying Talagrand's $\T_{2}$ inequality and not the logarithmic Sobolev inequality.
\end{abstract}
\section{Introduction}

Transport-entropy inequalities were introduced by Marton and Talagrand in the nineties \cite{Ma86,Ta96}. As their name indicates, this type of inequalities compare optimal transport costs in the sense of Monge-Kantorovich to the relative entropy functional (also called Kullback-Leibler divergence). Transport-entropy inequalities have deep connections to the concentration of measure phenomenon \cite{Ledoux-book,Go09}, to log-Sobolev type inequalities \cite{OV00,BGL01,GRS11bis}, or large deviation theory \cite{GL07,Go09}. They also directly appear in the definition proposed by Lott, Villani and Sturm of a metric measured space with positive Ricci curvature \cite{LV09,St06a}. The interested reader can consult the books \cite{Ledoux-book, Vi09} or the recent survey \cite{GL10} for an overview of their applications.

The purpose of this note is to give a necessary and sufficient condition for a large class of transport-entropy inequalities involving probability measures on the real line. 

Before presenting our main result, we first need to define transport costs and transport-entropy inequalities.
Let $\alpha:\R^+ \to\R^+$ be a cost function; the optimal transport cost between two probability measures $\mu,\nu$ is defined by
\begin{equation}\label{optimal transport}
\mathcal{T}_{\alpha}(\nu,\mu) =\inf \iint \alpha(|x-y|)\,\pi(dxdy),
\end{equation}
where the infimum runs over the set of couplings $\pi$ between $\mu$ and $\nu$, \textit{i.e} probability measures on $\R^2$ such that $\pi(dx\times \R)=\mu(dx)$ and $\pi(\R\times dy)=\nu(dy).$

A Borel probability measure $\mu$ on $\R$ is said to satisfy the transport-entropy inequality $\T_{\alpha}(a)$ for some $a>0$ if
$$\mathcal{T}_{\alpha(a\,\cdot\,)}(\nu,\mu)\leq \H(\nu\mid\mu),$$ for all  $\nu\in \mathcal{P}(\R)$ (the set of Borel probability measures on $\R$),
where $\alpha(a\,\cdot\,)$ denotes the cost function $t\mapsto \alpha(at)$ and where $\H(\nu\mid\mu)$ stands for the relative entropy of $\nu$ with respect to $\mu$. This latter is defined by $$\H(\nu\mid\mu)=\int \log \left( \frac{d\nu}{d\mu}\right)\,d\nu,$$ when $\nu$ is absolutely continuous with respect to $\mu$ and $\infty$ otherwise. For instance, the celebrated Talagrand's inequality $\T_2$ enters this family of inequalities. We recall that $\mu$ satisfies $\T_2(C)$ if
\begin{equation}\label{T_2}
\mathcal{T}_2(\nu,\mu)\leq C \H(\nu \mid \mu),\qquad \forall \nu\in \mathcal{P}(\R),
\end{equation}
where $\mathcal{T}_2$ is an abbreviated notation for $\mathcal{T}_{x^2}.$ With the definition introduced above, $\T_2(C)$ holds if and only if $\T_{x^2}\left(1/\sqrt{C}\right)$ holds.

In all the paper, we will use the following notation. The cumulative distribution function $F_{\nu}$ of a probability measure $\nu$ on $\R$ is the right continuous and non-decreasing function defined by $$F_{\nu}(x)=\nu(-\infty,x],\qquad \forall x\in \R.$$ The generalized inverse of $F_{\nu}$ is defined by  $$F_{\nu}^{-1}(u)=\inf\{x\in \R; F(x)\geq u\} \in \R\cup\{\pm\infty\},\qquad \forall u\in [0,1].$$ If $\mu$ is a probability measure with no atom and $\nu$ is another probability measure we will denote by $T_{\mu,\nu}$ the map defined by 
\begin{equation}\label{transport map}
T_{\mu,\nu}=F_{\nu}^{-1}\circ F_{\mu}.
\end{equation}
It is well known that $T_{\mu,\nu}$ is the only one non-decreasing and left-continuous function that pushes forward $\mu$ onto $\nu$, that is to say
$$\int f\,d\nu = \int f\circ T_{\mu,\nu}\,d\mu,$$
for all bounded measurable $f:\R\to\R.$\\
In what follows the exponential distribution 
\begin{equation}\label{exponential}
\mu_{1}(dx)=e^{-|x|}\,dx/2
\end{equation}
will play a central role. 

\pagebreak
In this paper, we will say that that a Borel probability $\mu$ on $\R$ satisfy Poincar\'e inequality with the constant $\lambda>0$ if
\begin{equation}\label{Poincare}
\lambda \mathrm{Var}_{\mu}(f) \leq \int |\nabla f|^2\,d\mu,\qquad \forall f\text{ Lipschitz},
\end{equation} 
where 
\begin{equation}\label{nabla}
|\nabla f|(x)= \limsup_{y\to x} \frac{|f(y)-f(x)|}{|y-x|},\qquad \forall x\in \R.
\end{equation}
Note that when $f$ is differentiable at $x$, then $|\nabla f|(x)=|f'(x)|$. Proposition \ref{equiv poinc} clarifies this definition of the Poincar\'e inequality.

The following theorem is our main result. It characterizes transport-entropy inequalities $\T_{\alpha}$ for convex functions $\alpha$ which are quadratic near $0$. 

\begin{thm}\label{equivalence 2}Let $\mu$ be a Borel probability measure on $\R$ and $\alpha:\R^+\to\R^+$ be a convex function such that $\alpha(t)=t^2$ for all $t\leq h$.
The following propositions are equivalent
\begin{enumerate}
\item There is some $a>0$ such that $\mu$ verifies $\T_{\alpha}(a)$.
\item There are $\lambda>0$ and $d>0$ such that\\ 
(i) $\mu$ verifies Poincar\'e inequality with constant $\lambda$ and\\
(ii) the map $T:=T_{\mu_{1},\mu}$ sending $\mu_{1}$ on $\mu$ verifies 
\begin{equation}\label{contraction}
|T(x)-T(y)|\leq \frac{1}{d} \alpha^{-1}(h^2+|x-y|),\qquad \forall x,y \in \R.
\end{equation}
\end{enumerate}
Moreover, there exist two positive constants $\kappa_{1},\kappa_{2}$ depending only on $h$ such that the optimal constants $a_{\mathrm{opt}}, \lambda_{\mathrm{opt}}, d_{\mathrm{opt}}$ are related as follows:
$$ \kappa_{1} \min\left(\sqrt{\lambda_{\mathrm{opt}}};d_{\mathrm{opt}}\right)\leq a_{\mathrm{opt}}\leq \kappa_{2} \min\left(\sqrt{\lambda_{\mathrm{opt}}};d_{\mathrm{opt}}\right).$$
\end{thm}

In other words, the transport-entropy inequality $\T_{\alpha}$ carries two different informations: the existence of a spectral gap and a quantitative information on the way the exponential distribution $\mu_{1}$ is deformed in order to produce $\mu$. Theorem \ref{equivalence 2} improves the results obtained by the author in a preceding work \cite{Go07}, where different necessary or sufficient conditions were investigated (see Section 4.1 for a discussion). Here, a true equivalence is obtained. 

It is well known that an absolutely continuous probability measure $\mu$ on $\R$ verifies Poincar\'e inequality if and only if the following holds
\begin{gather}\label{Muckenhoupt}
A^+:= \sup_{x\geq m} \mu [x,\infty) \int_{m}^x \frac{1}{p(t)}\,dt <\infty,\\
A^-:= \sup_{x\leq m} \mu(-\infty,x)  \int_x^{m} \frac{1}{p(t)}\,dt <\infty,\notag
\end{gather}
where $p$ denotes the density of $\mu$ with respect to the Lebesgue measure and $m$ is a median of $\mu.$ This result follows from a similar necessary and sufficient condition for weighted Hardy's inequalities due to Muckenhoupt \cite{Muc72} (extending previous works by Artola, Talenti \cite{Ta69} and Tomaselli \cite{Tom69}). Moreover, it can be shown (see e.g \cite{ane}) that the optimal constant $\lambda_{\mathrm{opt}}$ in Poincar\'e inequality \eqref{Poincare} verifies
$$\max(A^- ; A^+)\leq 1/\lambda_{\mathrm{opt}} \leq 4 \max(A^- ; A^+),$$
with possible cases of equalities see \cite{Mi08}.

To complete Theorem \ref{equivalence 2}, we shall give in Section 4 an easy to check sufficient condition for the contraction property \eqref{contraction} for absolutely continuous $\mu$ with smooth density. This condition deals with the asymptotic behavior of the logarithm of the density of $\mu$.

Theorem \ref{equivalence 2} is satisfactory from a theoretical point of view. Its conclusion is reminiscent of the characterizations of different functional inequalities on the line by Bobkov and Houdr\'e \cite{BH00,BH97} and Bobkov and G\"otze \cite{BG99bis}.
Theorem \ref{equivalence 2} is also a useful tool for constructing examples illustrating borderline situations. We will use it in the last section  to give a new example of a probability measure which verifies Talagrand's $\T_2$ inequality but not the logarithmic Sobolev inequality. Contrary to the previous example given by Cattiaux and Guillin in \cite{CG06}, the tail behavior of the probability exhibited in the present paper is exactly Gaussian. In the same section, we will answer a question raised by Cattiaux and Guillin in \cite{CG06} about the equivalence of Talagrand's inequality to Gaussian concentration and Poincar\'e inequality. We will use Theorem \ref{equivalence 2} again to give an appropriate counterexample.

One of the main ingredient in the proof of Theorem \ref{equivalence 2} is the fact that optimal transport has a very simple structure in dimension one. The following theorem is very classical and goes back to the works by Hoeffding, Fr\'echet and Dall'Aglio \cite{H40,Fr60,DA56}. A proof can be found in the books by Villani \cite{Vi03} or Rachev-Ruschendorf \cite{RR98}.

\begin{thm}\label{HFDA-intro}
Let $\alpha : \R^+\to\R^+$ be a convex function such that $\alpha(0)=0$ and suppose that $\mu\in \mathcal{P}(\R)$ has no atom, then for all probability measure $\nu \in \mathcal{P}(\R)$ such that $\iint \alpha(|x-y|)\,\mu(dx)\nu(dy)<\infty$,  the map $T_{\mu,\nu}$ defined by \eqref{transport map} realizes the optimal transport of $\mu$ onto $\nu$. In other words, the coupling $\pi(dxdy)=\delta_{T_{\mu,\nu}(x)}\mu(dx)$ achieves the infimum in \eqref{optimal transport} and so $$\mathcal{T}_\alpha(\nu,\mu)=\int \alpha(|x-T_{\mu,\nu}(x)|)\,\mu(dx).$$
\end{thm}

An immediate consequence of Theorem \ref{HFDA-intro} is that the optimal transport cost $\mathcal{T}_{\alpha}(\nu,\mu)$ is linear with respect to $\alpha$ on the convex cone of non-negative convex cost functions $\alpha$ vanishing at $0$: in particular, if $\alpha=\alpha_1+\alpha_{2}$ with $\alpha_{i}:\R^+\to\R^+$ a convex function, then
$$\mathcal{T}_{\alpha}(\nu,\mu)=\mathcal{T}_{\alpha_{1}}(\nu,\mu)+\mathcal{T}_{\alpha_{2}}(\nu,\mu).$$
This property is really specific of the dimension one. In general, one only has the trivial inequality 
\begin{center}$\mathcal{T}_\alpha\geq \mathcal{T}_{\alpha_1}+\mathcal{T}_{\alpha_2}.$\end{center}

To prove Theorem \ref{equivalence 2}, we shall use this observation with a decomposition of $\alpha$ into a function $\alpha_1$ which is quadratic near $0$ and then linear and a function $\alpha_2$ which vanishes in a neighborhood of $0$ and has the same growth as $\alpha$. The transport inequality $\T_\alpha$ is thus equivalent to the realization of both $\T_{\alpha_1}$ and $\T_{\alpha_2}$. The transport-entropy inequality $\T_{\alpha_{1}}$ is equivalent to Poincar\'e inequality as proved by Bobkov, Gentil and Ledoux \cite{BGL01} (see also Theorem \ref{BGL} below). We shall establish that $\T_{\alpha_2}$ is equivalent to the contraction condition \eqref{contraction}, which will complete the proof of Theorem \ref{equivalence 2}.

\bigskip

The paper is organized as follows. Section 2 is devoted to transport-entropy inequalities associated to functions $\alpha$ vanishing in a neighborhood of $0$. This class of transport-entropy inequalities have their own interest since they can be even be verified by discrete probability measures. We show the equivalence between these inequalities and contraction properties like \eqref{contraction}. In Section 3, we complete the proof of Theorem \ref{equivalence 2} following the strategy explained above. Section 4 is devoted to examples. The article ends with an appendix relating the definition we adopted of Poincar\'e inequality \eqref{Poincare} to other more classical formulations.

\section{Transport-entropy inequalities for costs vanishing in a neighborhood of $0$}
To begin with, let us observe that Talagrand's inequality $\T_{2}$ cannot be satisfied by a discrete probability measure of the form $$\mu= \sum_{k\in \N} \mu_{k} \delta_{k},$$
where the $\mu_{k}$'s are non-negative negative numbers of sum equal to $1$. Indeed, if a probability measure verifies $\T_{2}$ then it verifies Poincar\'e inequality \eqref{Poincare} (see for instance the proof of Theorem \ref{equivalence 2}), which excludes probabilities $\mu$ as above (unless it is a Dirac mass). 

In this section, we study transport-entropy inequalities associated to cost functions which are identically $0$ in a neighborhood of $0$. 
As we shall see, the interest of this type of cost functions is that the associated transport-entropy inequality can also be satisfied by discrete probability measures. Let us mention that inequalities of this type appeared also in a paper by Bonciocat and Sturm \cite{BS09} in their study of curvature of discrete metric spaces. 

In all what follows, $\beta:\R^+\to\R^+$ will be a convex function such that $\beta(t)=0$ for all $t\leq h,$ for some $h>0,$ and $\beta$ is increasing on $[h,\infty).$ The main result of this section is the following

\begin{thm}\label{equivalence 3}
A Borel probability measure $\mu$ on $\R$ verifies the transport-entropy inequality $\T_{\beta}(a)$ for some constant $a>0$ if and only if the transport map $T=T_{\mu_{1},\mu}$ sending the exponential distribution $\mu_{1}$ onto $\mu$ verifies the contraction property 
\begin{equation}\label{contraction bis}
|T(u)-T(v)|\leq \frac{1}{d}\beta^{-1}(|u-v|),\qquad \forall u\neq v\in \R.
\end{equation}
Moreover, the optimal constants $a_{\mathrm{opt}}$ and $d_{\mathrm{opt}}$ verify
$$ d_{\mathrm{opt}} \left(\frac{h}{9\beta^{-1}(2)}\right)\leq a_{\mathrm{opt}}\leq d_{\mathrm{opt}} \left(\frac{8\beta^{-1}(\log(3))}{h}\right).$$ 
\end{thm}
It is very easy to construct discrete probabilities enjoying a transport-entropy inequality $\T_{\beta}$. For example, consider the map $T:\R\to\N$ defined by $T(x)=\lceil \sqrt{x}\rceil$, for $x\geq0$ and $T(x)=0$ for $x\leq 0$, where $\lceil x\rceil$ is the smallest $k\in \N$ such that $x\leq k.$ It is clear that 
$$|T(x)-T(y)| \leq 1+\sqrt{|x-y|},\qquad \forall x, y\in \R.$$
Define $\mu$ as the image of $\mu_1$ under $T$. Since $T$ is left continuous, we have $T=T_{\mu_1,\mu}$ and so $\mu$ verifies the transport-entropy inequality $\T_{\beta_{2}}(a)$ for some constant $a$ with the cost function $\beta_{2}$ defined by $$\beta_{2}(x) = [x-1]_{+}^2,\qquad \forall x\geq 0.$$ 
In this example, $h=d_{\mathrm{opt}}=1$ and so the optimal constant $a_{\mathrm{opt}}$ verifies
$$\frac{1}{9(1+\sqrt{2})}\leq a_{\mathrm{opt}}\leq 8 (1+\sqrt{\log(3)}).$$

To prove Theorem \ref{equivalence}, we need to introduce some additional notation. 
Let $\mu$ be a probability measure on $\R$ which is not a Dirac mass  and define $$s_{\mu}=\inf\ \mathrm{Supp}(\mu)\qquad \text{and}\qquad t_{\mu}=\sup\ \mathrm{Supp}(\mu).$$ Let us define two families of probability measures $\{\mu_{x}^+\}$ and $\{\mu_{x}^-\}$ on $\R^+$ as follows:
$$\mu_{x}^+=\mathcal{L}(X-x|X> x),\quad \forall x<t_{\mu}$$
and
$$\mu_{x}^-=\mathcal{L}(x-X|X< x),\quad \forall x>s_{\mu},$$
where $X$ is a random variable with law $\mu.$\\In other words, for all bounded measurable function $f:\R\to\R$, 
$$\int f\,d\mu_{x}^+=\frac{\int_{(x,\infty)} f(u-x)\,\mu(du)}{\mu(x,\infty)}$$
and
$$\int f\,d\mu_{x}^-=\frac{\int_{(-\infty,x)} f(x-u)\,\mu(du)}{\mu(-\infty,x)}.$$
Define, for all $b\geq 0$
\begin{align}\label{K^+}
K^+(b)&=\sup_{t_{\mu}>x\geq m} \int_{0}^{\infty} e^{\beta(b u)}\,\mu_{x}^+(du)\in \R^+\cup\{\infty\}\\ K^-(b)&=\sup_{s_{\mu}<x\leq m} \int_{0}^{\infty} e^{\beta(b u)}\,\mu_{x}^-(du)\in \R^+\cup\{\infty\}\notag,
\end{align}
where $m$ is the median of $\mu$ defined by $m=F_{\mu}^{-1}(1/2),$ with the convention $\sup \emptyset =0.$

Theorem \ref{equivalence 3} follows immediately from the following improved version.

\begin{thm}\label{equivalence}
Let $\mu$ be a probability measure on $\R$ which is not a Dirac mass, and let $\mu_1$ be the two sided exponential distribution defined by \eqref{exponential}.
The following propositions are equivalent
\begin{enumerate}
\item There is $a>0$ such that $\mu$ verifies the transport inequality $\T_{\beta}(a).$
\item There are $b>0$ and $K>0$ such that $\max(K^-(b);K^{+}(b))\leq K$.
\item There is $c>0$ such that the map $S:\R \times [0,1] \to \R\cup \{\pm \infty\}$ defined by 
$$S(x,u)=F_{\mu_{1}}^{-1}\left( \mu (-\infty,x) + \mu\left(\left\{x\right\} \right)u\right),\qquad \forall x\in \R, \quad \forall u\in [0,1],$$
verifies
$$|S(x,u)-S(y,v)|\geq \beta(c|x-y|),\qquad \forall x,y\in \R,\quad \forall u,v\in [0,1].$$
\item There is $d>0$ such that the map $T:=T_{\mu_{1},\mu}$ defined by \eqref{transport map} which sends $\mu_{1}$ onto $\mu$ verifies
$$|T(u)-T(v)|\leq  \frac{1}{d}\beta^{-1}(|u-v|),\qquad \forall u\neq v\in \R.$$
\end{enumerate}
The constants are related in the following way:
\begin{center}
\begin{tabular}{lll}
(1) $\Rightarrow$ (2) \quad with \quad $b=a/2$ and $K= 3$.\\
(2) $\Rightarrow$ (3) \quad with \quad $c=b\left(\frac{h}{4\beta^{-1}(k)}\right)$ and $k=\log K.$\\
(3) $\Rightarrow$ (4) \quad with \quad $d=c$.\\
(4) $\Rightarrow$ (1) \quad with \quad $a=d\left(\frac{h}{9\beta^{-1}(2)}\right)$.
\end{tabular}
\end{center}
\end{thm}

Let us give an interpretation of the map $S$ appearing in condition (3). 
More generally, if $\mu$ and $\nu$ are arbitrary Borel probability measures on $\R$, we define the map $S_{\mu,\nu}:\R\times [0,1] \to \R\cup\{\pm \infty\}$ as follows:
\begin{equation}\label{transport map 2}
S_{\mu,\nu}(x,u) = F_{\nu}^{-1}\left( \mu(-\infty,x) + \mu\left(\left\{x\right\} \right)u\right),\qquad \forall x\in \R, \quad \forall u\in [0,1].
\end{equation}
Remark that in case $\mu$ has no atom, $S_{\mu,\nu}$ coincides with $T_{\mu,\nu}$ defined by \eqref{transport map}. As the following theorem explains, this map realizes the optimal transport of $\mu$ onto $\nu$.

\begin{thm}\label{HFDA}
Let  $\alpha:\R^+\to\R^+$ a convex cost function such that $\alpha(0)=0$ and $\mu,\nu$ be two probability measures on $\R$ such that $\int \alpha(|x-y|)\,\mu(dx)\nu(dy)<\infty$; then the coupling $\pi_{o}\in \mathcal{P}(\R^2)$ whose distribution function is given by
$$\pi_{o}((-\infty,x]\times (-\infty,y])=\min(F_{\mu}(x),F_{\nu}(y)),\qquad \forall x,y\in \R$$
achieves the infimum in the definition of $\mathcal{T}_{\alpha}(\nu,\mu).$
Moreover, if $X$ is a random variable with law $\mu$ and $U$ a random variable uniformly distributed on $[0,1]$ and independent of $X$, then
$$\pi_{o} = \mathrm{Law} (X,S_{\mu,\nu}(X,U)).$$
\end{thm}
Theorem \ref{HFDA} generalizes Theorem \ref{HFDA-intro}; we state it for completeness but it will not be used in the sequel. Note that the coupling $\pi_o$ remains optimal for a more general class of transport costs \cite{CSS76}.

During the proof of Theorem \ref{equivalence}, we will use the following simple technical lemma twice.
\begin{lem}\label{petit lemme}
Let $\beta:\R^+ \to\R^+$ be a convex function such that $\beta=0$ on $[0,h]$ and $\beta$ is increasing on $[h,\infty).$ Then, for all $b>0$ and $k>0$
$$[\beta(bv)-k]_+ \geq \beta (cv),\qquad \forall v\geq0,$$
where $c=b\left(\frac{h}{2\beta^{-1}(k)}\right).$
\end{lem}

\proof
If $v\leq h/c$, there is nothing to prove. If $v>h/c$, then since $\frac{b}{2} = c\frac{\beta^{-1}(k)}{h}> c,$ it holds
$$\beta(bv)\geq \beta(bv/2) + \beta(bv/2)\geq \beta(cv)+\beta(c\beta^{-1}(k)v/h) \geq \beta(cv)+k,$$
which proves the claim.
\endproof

\proof[Proof of Theorem \ref{equivalence}.] \ \\
\noindent(1) $\Rightarrow$ (2). According \cite[Proposition 8.3]{GL10} or \cite[Proposition 4.13]{Go10}, the assumed transport inequality implies the following inf-convolution inequality
$$\int e^{Qf}\,d\mu\int e^{-f}\,d\mu\leq 1,$$
for all function $f$ bounded from below, where 
$$Qf(y)=\inf_{z\in \R}\left\{f(z)+2\beta(a|y-z|/2)\right\}.$$
Consider the function $f_{x}$ which is $0$ on $(-\infty,x]$ and $\infty$ otherwise, then $Qf(y)=2\inf_{z\leq x} \beta(a|y-z|/2)$ and so
$Qf=0$ on $(-\infty,x]$ and $Qf(y)=2\beta(a(y-x)/2)$ on $(x,\infty).$ Applying the inequality above to $f_{x}$ thus yields
$$\left(\mu(-\infty,x] + \int_{(x,\infty)} e^{2\beta(a(y-x)/2)}\,\mu(dy)\right) \mu(-\infty,x]\leq 1.$$
From this follows that if $x\geq m$
$$\int_{0}^{\infty}e^{2\beta(au/2)}\,\mu_{x}^+(du)\leq\frac{1}{\mu(-\infty,x]}+1\leq 3.$$
So $K^+(a/2)\leq 3$ and similarly $K^-(a/2)\leq3.$\\

\noindent(2) $\Rightarrow$ (3). To prove (3) we can first restrict to the case $y>x$ and then using the monotonicity of $S$ we can further assume that $v=0$ and $u=1$. So Property (3) is equivalent to the following one
\begin{equation}\label{3 bis}
S(y,0)-S(x,1)\geq \beta(c(y-x)),\qquad \forall y>x.
\end{equation}
To establish \eqref{3 bis} it is enough to consider the cases $y>x\geq m$ and $m\geq y>x$. Namely, suppose that \eqref{3 bis} is true with a constant $\tilde{c}$ for these two particular cases, and consider $y>m>x$. Then, it holds
\begin{align*}
S(y,0)-S(x,1)&\geq S(y,0)-S(m,1)+S(m,0)-S(x,1)\\
& \geq \beta(\tilde{c}(y-m))+\beta(\tilde{c}(m-x))\\
&\geq \beta(\tilde{c}(y-x)/2),
\end{align*}
where the first inequality comes from the fact that $S(m,1)-S(m,0)\geq 0$ and the last one from the monotonicity of $\beta\geq 0$ and the inequality $\min(m-x ; y-m)\geq (y-x)/2.$\\
Let us check \eqref{3 bis} when $m\leq x<y$. Since $F_{\mu_{1}}^{-1}(t)=-\log(2(1-t))$ for $t\in [1/2 ; 1)$ and $F_{\mu_{1}}^{-1}(t)=\log(2t)$ for $t\in (0,1/2)$, it holds (since $\mu (-\infty,x]\geq1/2$)
\begin{align*}
S(y,0)-S(x,1)&= -\log(2(1-\mu(-\infty,y))+\log(2(1-\mu(-\infty,x])) \\
&= -\log( \mu_{x}^+ [y-x,\infty)).
\end{align*}
Since $K^+(b)\leq K$, Markov inequality implies that $\mu_{x}^+ ([u,\infty)) \leq Ke^{-\beta(bu)}$ for all $u>0.$ So, 
$$S(y,0)-S(x,1)\geq [\beta(b(y-x)) - \log(K)]_{+} \geq \beta(\tilde{c}(y-x)),\qquad \forall y>x\geq m,$$
with $\tilde{c}=b\left(\frac{h}{2\beta^{-1}(k)}\right)$ and $k=\log(K)$, where the second inequality follows from Lemma \ref{petit lemme}.
Reasoning exactly as above we show that the same inequality holds when $x<y\leq m$. So, according to what precedes, (3) holds with $c=\tilde{c}/2.$\\

\noindent(3) $\Rightarrow$ (4). By assumption, it holds 
$$|F_{\mu_{1}}^{-1} (\mu(-\infty,x) +\mu(\{x\}) u)- F_{\mu_{1}}^{-1} (\mu(-\infty,y) +\mu(\{y\}) v)|\geq \beta(c|x-y|),$$
for all $x,y\in \R$ and $u,v\in [0,1].$ Let us apply this inequality to $x=F_{\mu}^{-1}(s)$ and $y=F_{\mu}^{-1}(t)$ with $s,t \in (0,1).$
It is easy to check that $$\mu(-\infty, F_{\mu}^{-1}(s))\leq s \leq \mu(-\infty, F_{\mu}^{-1}(s)].$$
So choosing properly $u$ and $v$ yields
$$|F_{\mu_{1}}^{-1} (s) - F_{\mu_{1}}^{-1} (t)| \geq \beta\left(c|F_{\mu}^{-1}(s)-F_{\mu}^{-1}(t)|\right).$$
Finally applying this inequality to $s=F_{\mu_{1}}(z)$ and $t=F_{\mu_{1}}(w)$ gives the desired inequality.\\

\noindent(4) $\Rightarrow$ (1). According to \cite{Ma91}, the exponential distribution $\mu_{1}$ verifies the following inf-convolution inequality
\begin{equation}\label{Maurey}
\int e^{\overline{Q}g}\,d\mu_1\leq e^{\int g\,d\mu_1},
\end{equation}
for all bounded measurable $g$, where $\overline{Q}g(x)=\inf_{y\in \R}\left\{g(y)+ \beta_1(|x-y|)\right\},$
with $\beta_1(x)=\frac{1}{36} x^2$ if $0\leq x\leq 4$ and $\beta_1(x)=\frac{2}{9}(x-2)$ if $x\geq 4$. In particular, $\beta_1(x)\geq \frac{2}{9}[x-2]_+$, for all $x\geq0.$

According to (4), 
$$\frac{2}{9}\left[\beta(d|T(x)-T(y)|)-2\right]_+\leq \beta_1(|x-y|),\qquad \forall x,y\in \R.$$
According to Lemma \ref{petit lemme},  $\frac{2}{9}\left[\beta(dv)-2\right]_+ \geq \beta(av)$, with 
$a=d\left(\frac{h}{9\beta^{-1}(2)}\right).$

So, defining $$Qf(x)=\inf_{y\in \R}\left\{f(y)+\beta(a|x-y|)\right\},$$
we have
$$\left(Qf\right)(T(x))\leq \inf_{z\in \R}\left\{f(T(z))+\beta(a|T(x)-T(z)|)\right\}\leq \overline{Q}(f\circ T).$$
So applying \eqref{Maurey} to $g=f\circ T$, we get
$$\int e^{Qf}\,d\mu\leq e^{\int f\,d\mu},$$
for all bounded measurable $f$. According to Bobkov and G\"otze dual characterization \cite{BG99} (see also \cite{GL10}), we conclude that $\mu$ verifies $\T_\beta(a).$
\endproof

\section{Proof of Theorem \ref{equivalence 2}}
According to Bobkov, Gentil and Ledoux \cite{BGL01}, the Poincar\'e inequality is equivalent to a family of transport-entropy inequalities involving the cost functions $\alpha_{1}^h$ defined by 
$$\alpha^h_{1}(t)=t^2,\quad\text{ if } 0\leq t\leq h\quad\quad\text{and}\quad\quad \alpha^h_{1}(t)=2ht-h^2,\quad\text{ if } t\geq h.$$
More precisely, 
\begin{thm}[Bobkov-Gentil-Ledoux \cite{BGL01}]\label{BGL}
Let $\mu$ be a Borel probability measure on $\R$. The following propositions are equivalent:
\begin{enumerate}
\item There is $\lambda>0$ such that $\mu$ verifies the Poincar\'e inequality \eqref{Poincare} with the constant $\lambda.$
\item There are $a,h>0$ such that $\mu$ verifies $\T_{\alpha_1^{h}}(a).$
\end{enumerate}
The constants are related as follows\\
$(1)\Rightarrow (2)$ with $a=\frac{1}{2\sqrt{K(c)}}$, and $h=c\sqrt{K(c)}$ where $K(c)=\frac{1}{2\lambda}\left(\frac{2\sqrt{\lambda}+c}{2\sqrt{\lambda}-c}\right)^2e^{c\sqrt{5/\lambda}}$, for all $c\in [0,2\sqrt{\lambda}).$\\
$(2)\Rightarrow (1)$ with $\lambda=2a^2.$
\end{thm}
The preceding theorem is stated in dimension one only, but it is true in any dimension. 
\proof
The implication $(2)\Rightarrow (1)$ is true on any metric space (see \cite{GRS12bis}).
We refer to \cite{BGL01} or \cite{Vi09} for the proof of $(1)\Rightarrow (2)$ in the case when $\mu$ is absolutely continuous with respect to Lebesgue. In what follows, we show that this implication is still true when $\mu$ is not.

Let $\mu$ be a Borel probability measure on $\R$. For all $\sigma>0$, let $\gamma_\sigma = \mathcal{N}(0,\sigma^2)$ be a centered Gaussian distribution with variance $\sigma^2$ and define $\mu_\sigma=\mu\ast\gamma_\sigma$. The probability $\gamma_\sigma$ verifies the Poincar\'e inequality with the constant $1/\sigma^2.$ According to the well known tensorization property of Poincar\'e inequality \cite{Ledoux-book}, it is not difficult to check that the product measure $\mu\otimes \gamma_\sigma$ verifies the following inequality
$$\mathrm{Var}_{\mu\otimes \gamma_\sigma} (g) \leq  \int \frac{1}{\lambda^2}|\nabla_x g|^2(x,y) + \sigma^2|\nabla_y g|^2(x,y)\,\mu\otimes \gamma_\sigma (dxdy),$$
for all Lipschitz function $g:\R^2\to\R$.
Considering functions $g$ of the form $g(x,y)=f(x+y)$, we obtain
$$\mathrm{Var}_{\mu_\sigma} (f) \leq \left(\frac{1}{\lambda^2} + \sigma^2\right)\int |\nabla f|^2(z)\,\mu_\sigma (dz).$$
So $\mu_\sigma$ verifies Poincar\'e with the constant $\lambda_\sigma=\left(\frac{1}{\lambda^2} + \sigma^2\right)^{-1}$. Since $\mu_\sigma$ is absolutely continuous, we can conclude applying  \cite{BGL01} that $\mu_\sigma$ verifies the family of transport-entropy inequalities $\T_{\alpha_1^h}(a)$ with $a,h$ satisfying the constraints given in Theorem \ref{BGL}. Since $\mu_\sigma \to \mu$ for the weak topology and $\lambda_\sigma \to \lambda$, when $\sigma$ goes to $0$, it is not difficult to see that $\mu$ verifies the transport-entropy inequalities $\T_{\alpha_1^h}(a)$ for $a$ and $h$ in the good range. (This last step is easier to check on the dual form of Bobkov-G\"otze.)
\endproof

We are now ready to prove Theorem \ref{equivalence 2} using the decomposition trick explained in the introduction.
\proof[Proof of Theorem \ref{equivalence 2}]
(1) $\Rightarrow$ (2). Observe that $\alpha(a\,\cdot\,) \geq \alpha^h_{1}(a\,\cdot\,)$. This inequality is immediate when $t\leq 1/a$ and results from the convexity of $\alpha$ when $t\geq1/a$. Therefore, $\mu$ verifies $\T_{\alpha^h_{1}}(a)$ and so the Poincar\'e inequality with the constant $2a^2.$ On the other hand, the inequality $\alpha(a\,\cdot\,)\geq \alpha_2(a\,\cdot\,)$, with $\alpha_2=[\alpha-h^2]_{+}$ implies that $\mu$ verifies $\T_{\alpha_2}(a)$. According to Theorem \ref{equivalence}, we conclude that the transport map $T$ enjoys \eqref{contraction} with $d= a\left(\frac{h}{8\alpha^{-1}(h^2+\log(3))}\right).$ So
\begin{equation}\label{majoration}
a_{\mathrm{opt}}\leq \min\left(d_{\mathrm{opt}};\sqrt{\lambda_{\mathrm{opt}}}\right) \max \left(\frac{1}{\sqrt{2}} ; \frac{8\alpha^{-1}(h^2+\log(3))}{h}\right)
\end{equation}

(2) $\Rightarrow$ (1). Let $c_{o}$ be such that $c_{o}\sqrt{K(c_{o})}=h$, then, according to Theorem \ref{BGL}, $\mu$ verifies $\T_{\alpha_{1}^h}(a_{1})$ with $a_{1}=\frac{1}{2\sqrt{K(c_{o})}}=\frac{c_{o}}{2h}$. It is not difficult to check that $a_{1}\geq \sqrt{\lambda} \frac{\kappa}{1+\kappa h}$, with $\kappa=\left(\sqrt{2}e^{-\sqrt{5}}\right)/4.$ Define $\alpha_2=[\alpha -h^2]_{+}$; according to Theorem \ref{equivalence}, $\mu$ verifies $\T_{\alpha_2}(a_{2})$, with $a_{2}=d\left(\frac{h}{9\alpha^{-1}(h^2+2)}\right)$. Observe that the function $\alpha-\alpha_1^h:\R^+\to\R^+$ is convex and verifies the inequality $\alpha-\alpha_{1}^h\leq \alpha_2.$ This inequality is clear on $[0,h]$ and for $t\geq h$, it holds
$$\alpha(t)-\alpha^h_{1}(t) = \alpha(t)-2ht +h^2 = \alpha_2(t)-2h(t-h)\leq \alpha_2(t).$$
So defining $a=\min(a_1;a_2)$ and applying Theorem \ref{HFDA-intro}, we thus have
$$\mathcal{T}_{\alpha(a\,\cdot\,)}(\nu,\mu)\leq \mathcal{T}_{\alpha_1^h(a_1\,\cdot\,)}(\nu,\mu)+\mathcal{T}_{\alpha_2(a_2\,\cdot\,)}(\nu,\mu),\qquad \forall \nu\in \mathcal{P}(\R).$$
So, $\mu$ verifies $\T_{\alpha / 2}(a)$ wich is stronger than $\T_\alpha(a/2)$.
So 
\begin{equation}\label{minoration}
a_{\mathrm{opt}} \geq \frac{1}{2}\min\left(d_{\mathrm{opt}};\sqrt{\lambda_{\mathrm{opt}}}\right) \min\left(\frac{\kappa}{1+\kappa h} ; \frac{h}{9\alpha^{-1}(h^2+2)}\right).
\end{equation}
To complete the proof, observe that since $\alpha$ is convex and $\alpha(x)=x^2$ on $[0,h]$, one has
$\alpha(x)\geq 2xh - h^2$ for all $x\geq0$. Therefore, $\alpha^{-1} (y) \leq \frac{y+h^2}{2h}$, for all $y\geq 0.$ Plugging this upper bound into \eqref{majoration} and \eqref{minoration} yields
$$\kappa_{1} \min\left(d_{\mathrm{opt}};\sqrt{\lambda_{\mathrm{opt}}}\right)\leq a_{\mathrm{opt}} \leq \kappa_{2} \min\left(d_{\mathrm{opt}};\sqrt{\lambda_{\mathrm{opt}}}\right),$$
with 
$$\kappa_{1}=\frac{1}{2}\min\left(\frac{\kappa}{1+\kappa h} ; \frac{h^2}{9(h^2+1)}\right),$$
and $$\kappa_{2}= \frac{4}{h^2} (2h^2 +\log(3)).$$
\endproof
\section{Examples}
This section is devoted to examples. First we recall the result obtained in \cite{Go07} and make the link with the present paper.  After that, we give a general sufficient condition for transport-entropy inequalities which holds for absolutely continuous distributions with smooth densities. We end the section by showing how Theorem \ref{equivalence 2} can be used to construct borderline examples, typically a probability enjoying $\T_{2}$ but not the logarithmic Sobolev inequality.

\subsection{Connection with \cite{Go07}}
Let us make the connection between \cite{Go07} and the present paper. Let us recall that a probability measure $\mu$ on $\R$ verifies Cheeger's inequality, if
$$\int |f-m|\,d\mu \leq \int |\nabla f|\,d\mu,\qquad \text{for all $f$ Lipschitz,}$$
where $m$ is a median of $\mu$ and $|\nabla f|$ is defined by \eqref{nabla}. Cheeger's inequality is known to be strictly stronger than Poincar\'e inequality. For probability distributions on $\R$, it was proved by Bobkov and Houdr\'e \cite{BH97} that Cheeger's inequality holds if and only if the transport map $T_{\mu_{1},\mu}$ is Lipschitz.

In \cite{Go07}, we obtained the following incomplete characterization
\begin{thm}\label{Thm Go07}Let $\mu$ be an absolutely continuous distribution on $\R$ verifying Cheeger's inequality; $\mu$ verifies $\T_{\alpha}(a)$ for some $a>0$ if and only if there is some $b>0$ such that $\max(K^-(b);K^+(b))<\infty$, where $K^\pm$ are defined by \eqref{K^+} (with $\beta=\alpha$).
\end{thm}
It is not difficult to construct a probability verifying for example $\T_{2}$ and not Cheeger's inequality (and thus which is not covered by Theorem \ref{Thm Go07}).
For example, consider the probability $\nu(dx)=\frac{1}{Z} |x|^re^{-|x|}\,dx$ for some $r\in (0,1)$. One can check that $\nu$ verifies Muckenhoupt's conditions \eqref{Muckenhoupt} and so Poincar\'e. Let $T_{1}$ be the transport map $T_{\mu_{1},\nu}$. Writing $F_{\nu}(x)=F_{\mu_{1}}(T_{1}^{-1}(x))$ and taking the derivative at $x=0$, we see that $T_{1}'(x)\to\infty$ when $x\to 0$ and so $T_{1}$ is not Lipschitz. According to Bobkov and Houdr\'e \cite{BH97}, it follows that $\nu$ does not verify Cheeger's inequality (this example is taken from \cite{BH97}). 
Now, consider $T_{2}(x)=\mathrm{sign}(x)\min(|x| ; \sqrt{|x|})$ and define $\mu$ as the image of $\nu$ under $T_{2}.$ 
We claim that $\mu$ verifies Talagrand's inequality $\T_{2}$ and not Cheeger's inequality.
Indeed, since $\nu$ verifies Poincar\'e inequality, one concludes from Theorems \ref{BGL} and \ref{equivalence 2} that 
$$|T_{1}(x)-T_{1}(y)|\leq a + b |x-y|,\qquad \forall x,y \in \R$$
On the other hand,
$$|T_{2}(x)-T_{2}(y)|\leq 2\sqrt{|x-y|},\qquad \forall x,y\in \R.$$
Combining these two inequalities we see that $T=T_{2}\circ T_{1}$ verifies
$$|T(x)-T(y)| \leq 2\sqrt{a+b|x-y|},\qquad \forall x,y \in \R.$$
Moreover, since $T_{2}$ is $1$-Lipschitz, $\mu$ verifies Poincaré inequality and so according to Theorem \ref{equivalence 2} $\mu$ verifies $\T_{2}$. Finally,
$T'(x)=T_{2}'(T_{1}(x))T_{1}'(x)\to \infty$ when $x\to 0$ and so $\mu$ does not verify Cheeger's inequality.
 
\subsection{A general criterion on the density.} We recall below a sufficient condition obtained by the author in \cite{Go07} that ensures that a probability on $\R$ with a smooth density verifies a transport-entropy inequality.

\begin{thm}\label{examples}Suppose that $\alpha:\R^+ \to \R^+$ is a convex function of class $\mathcal{C}^2$ such that $\alpha(t)=t^2$ for small values of $t$ and verifying the following regularity assumption: $\frac{\alpha''(t)}{(\alpha'(t))^2}\to0$ when $t\to\infty.$ Let $\mu$ be an absolutely continuous probability measure on $\R$ with a density of the form $d\mu(x)=e^{-V(x)}\,dx$, where $V:\R\to\R$ is a function of class $\mathcal{C}^2$ such that $\frac{V''(t)}{(V'(t))^2}\to 0$ as $t \to\infty.$ If $V$ is such that there is $\lambda>0$ such that
\begin{equation}\label{liminf}
\liminf_{x\to\pm \infty} \frac{|V'(x+m)|}{\alpha'(\lambda |x|)}>0,
\end{equation}
where $m$ is the median of $\mu$, then $\mu$ verifies the transport-entropy inequality $\T_{\alpha}(a)$ for some $a>0.$
\end{thm}
Note that in the quadratic case, condition \eqref{liminf} was first obtained by Cattiaux and Guillin in \cite{CG06}.
The proof of \cite{Go07} goes as follows: using a classical asymptotic analysis, we show that the condition \eqref{liminf} ensures that $\max(K^-(b);K^+(b))$ is finite for $b$ small enough. On the other hand, the condition $\liminf_{x\to \pm \infty} |V'(x)|>0$ (which is implied by \eqref{liminf}) is enough to have Cheeger's inequality. The conclusion follows from Theorem \ref{Thm Go07}.

Let us mention that multidimensional generalizations of condition \eqref{liminf} were proposed in \cite{Go10} or in \cite{CGW10}.
In the one dimensional case, we do not know if it is possible to use Theorem \ref{equivalence 2} to significantly enlarge the class of examples given in Theorem \ref{examples}. 

\subsection{Counterexamples.} Our main result Theorem \ref{equivalence 2} enables us to exhibit new examples of probability measures clarifying the links between Talagrand's inequality \eqref{T_2} and the logarithmic Sobolev inequality, the Poincar\'e inequality \eqref{Poincare} and Gaussian concentration.

Let us recall that a Borel probability measure $\mu$ on $\R$ is said to verify the logarithmic Sobolev inequality if 
\begin{equation}\label{LSI}
\mathrm{Ent}_\mu(f^2):= \int f^2\log\left(\frac{f^2}{\int f^2\,d\mu}\right)\,d\mu \leq C \int |\nabla f|^2\,d\mu,
\end{equation}
for all $f$ Lipschitz, with $|\nabla f|$ defined by \eqref{nabla}. The known hierarchy between the above mentioned inequalities is the following:
$$\text{Log-Sobolev}\quad\Rightarrow\quad \T_2\quad \Rightarrow\quad \text{Poincar\'e}.$$
This chain of implications was first established by Otto and Villani in \cite{OV00} on Riemannian manifolds (see also \cite{BGL01}); it is true in a general framework \cite{GRS12}.

\subsubsection{A probability measure verifying $\T_2$ and not the logarithmic Sobolev inequality.}
In \cite{CG06}, Cattiaux and Guillin were the first to show that Talagrand's inequality was not equivalent to Log-Sobolev. They proved that the probability measure $\mu_{CG}$ defined on $\R$ by
$$\mu_{CG}(dx)= \frac{1}{Z} \exp( -|x|^3 - |x|^\beta -3x^2\sin^2(x))\,dx,\qquad \text{with}\quad 2<\beta<5/2,$$
verifies $\T_2$ but not the logarithmic Sobolev inequality. Our purpose is to produce another example whose tail distribution is exactly Gaussian.

Let us define a probability measure $\mu$ on $\R$ as the image of the exponential distribution $\mu_1(dx)=\exp(-|x|)\,dx/2$ under the map $T:\R\to\R$ defined as follows:  $T$ is odd, continuous and for all $k\in \N$, $T(x)= k$ on the interval $[k^2,(k+1)^2-1]$ and affine on $[(k+1)^2-1,(k+1)^2].$ 
We claim that this probability $\mu$ do the job. First, observe that $\mu$ verifies Poincar\'e inequality. This follows immediately from the fact that $T$ is $1$-Lipschitz. Moreover, it easily follows from the definition of $T$ that the following inequality holds: for all $y\geq x \in \R$
\begin{align*}
T(y)-T(x)&\leq 2+\mathrm{Card} \{ k\in \mathbb{N};  k^2 \in [x,y]\}+ \mathrm{Card} \{ k\in \mathbb{N};  -k^2 \in [x,y]\}\\
&\leq 4+2\sqrt{y-x}.
\end{align*}
According to Theorem \ref{equivalence 2}, we conclude that $\mu$ verifies $\T_2$. (Note that Theorem \ref{equivalence 2} actually applies because $T=F_{\mu}^{-1}\circ F_{\mu_1}$.)

To show that $\mu$ does not verify the logarithmic Sobolev inequality, we shall use the following criterion due to Bobkov and G\"otze \cite{BG99} (see also \cite{BR03}):
\begin{thm}  Let $\mu$ be a Borel probability measure on $\R$ and let $p:\R\to \R^+$ be the density of the absolutely continuous part of $\mu$. The probability $\mu$ verifies the logarithmic Sobolev inequality
\begin{equation}\label{LSI bis}
\mathrm{Ent}_\mu (f^2)\leq C \int |f'|^2(x) p(x)\,dx,\qquad \forall f \text{ Lipschitz,}
\end{equation}
if and only if 
$$D^+=\sup_{x\geq m} \mu[x,\infty)\log\left(\frac{1}{\mu[x,\infty)}\right) \int _m ^x \frac{1}{p(t)}\,dt <\infty$$
and 
$$D^-=\sup_{x\leq m} \mu(-\infty,x)\log\left(\frac{1}{\mu(-\infty,x)}\right) \int _x^m \frac{1}{p(t)}\,dt <\infty,$$
where $m$ is any median of $\mu.$ Moreover, the optimal constant $C_{LS}$ in \eqref{LSI} is such that $$c_1 \max(D^-;D^+)\leq C_{LS}\leq c_2 \max(D^-;D^+),$$
where $c_1,c_2$ are universal constants.
\end{thm}
\begin{rem} We refer to Proposition \ref{equiv poinc} for the relation between \eqref{LSI} and \eqref{LSI bis}. In particular, the probability $\mu$ defined above enters the class of probability measures for which \eqref{LSI} and \eqref{LSI bis} are equivalent.
\end{rem}
Let us come back to our example and show that the probability $\mu$ constructed above does not verify the logarithmic Sobolev inequality. We will show that $D^+=\infty.$
Let $f:\R^+\to\R$ be a bounded measurable function; then it holds
\begin{align*}
\int_0^{\infty} f\,d\mu &= \frac{1}{2}\int_{0}^{\infty} f\circ T(x)e^{-x}\,dx\\
& = \frac{1}{2}\sum_{k=1}^{\infty} f(k) \left(e^{-k^2}-e^{-(k+1)^2+1}\right) + \frac{1}{2}\sum_{k=0}^{\infty} \int_{(k+1)^2-1}^{(k+1)^2}f(u+k+1-(k+1)^2)e^{-u}\,du\\
&= \frac{1}{2}\sum_{k=1}^{\infty} f(k) \left(e^{-k^2}-e^{-(k+1)^2+1}\right) + \int_{0}^{\infty} f(t)p(t)\,dt,\end{align*}
where 
$$p(t)=\frac{e^{-t}}{2} \sum_{k=0}^{\infty} \mathbf{1}_{(k,k+1)}(t)e^{(k+1)-(k+1)^2},\qquad \forall t\geq0.$$
is the density of the absolutely continuous part of $\mu$ on $\R^+$. Observe also that the median of $\mu$ is $0$ and that for all $n\in \N$, $\mu[n,\infty) = \mu_1[n^2,\infty)=\frac{1}{2}e^{-n^2}.$ After some calculations, we get
\begin{align*}
D_n^+:&=\mu[n,\infty) \log\left(\frac{1}{\mu[n,\infty)} \right)\int_0^n\frac{1}{p(t)}\,dt\\
&=e^{-n^2}(n^2+\log(2))(1-1/e)\sum_{k=1}^{n}e^{k^2}\\
&\geq e^{-n^2}(n^2+\log(2))(1-1/e) \int_0^n e^{t^2}\,dt.
\end{align*}
Observing that,
$$\int_0^n e^{t^2}\,dt \geq \int_0^n \frac{t}{n}e^{t^2}\,dt= \frac{1}{2n}(e^{n^2}-1),$$
we conclude that $D_n^+ \to \infty$, when $n\to\infty$, and so $D^+=\infty.$ This completes the proof that $\mu$ does not verify the logarithmic Sobolev inequality. 

\begin{rem}If one wants to construct a counterexample $\tilde{\mu}$ absolutely continuous with respect to the Lebesgue measure, it suffices to replace in the definition of $T$ the constant steps by linear steps with small slope.
\end{rem}

\subsubsection{A probability with a Gaussian tail verifying Poincar\'e inequality and not $\T_2$.}
To motivate the construction of this probability, let us say a word on the tightening of functional inequalities. Recall that an absolutely continuous probability measure $\mu$ on $\R^n$ verifies the defective logarithmic Sobolev inequality if there are some constants $C,D\geq0$ such that
$$\mathrm{Ent}_\mu(f^2)\leq C \int |\nabla f|^2\,d\mu + D \int f^2\,d\mu,$$
for all $f:\R^n\to\R$ Lipschitz. A very classical result states that if $\mu$ verifies a defective logarithmic Sobolev inequality with constants $C,D$ and a Poincar\'e inequality with constant $\lambda$, then it verifies the logarithmic Sobolev inequality with a constant that can be expressed in terms of $C,D$ and $\lambda$.  Up to a subtle centering argument due to Rothaus \cite{Ro85}, this tightening result is intuitively clear. The tightening recipe 
\begin{center}
``defective functional inequality + Poincar\'e inequality = tight functional inequality'' 
\end{center}
appears to be very general, and holds for a large class of functional inequalities (see e.g \cite{BCR06,BK08}). A natural question is to ask if this tightening principle holds for transport-entropy inequalities. 

Let us say that a probability $\mu$ on $\R^n$ equipped with its standard Euclidean norm $\|\,\cdot\,\|_2$ verifies the defective transport-entropy inequality $\T_2$ if there are $C,D\geq 0$ such that 
$$\mathcal{T}_2(\nu,\mu)\leq C \H(\nu\mid\mu) +D,$$
for all probability measure $\nu.$ (The transport cost $\mathcal{T}_2(\nu,\mu)$ is defined as the infimum of $\mathbb{E}\left[ \|X-Y\|_2^2\right]$ over all the possible random variables $X,Y$ with respective law $\mu$ and $\nu$.) This defective $\T_2$ inequality has been characterized in various places (\cite{CG06,BV05,Go06}). It has been shown that it was equivalent to Gaussian concentration or equivalently to the finiteness of $\int e^{\varepsilon \|x\|_{2}^2}\,\mu(dx)$ for some $\varepsilon>0.$ Therefore, if the tightening principle was true for transport-entropy inequalities, then we would have the following equation
\begin{equation}\label{conjecture}
\T_2 \quad =\quad  \text{Poincar\'e} \quad +\quad \int e^{\varepsilon \|x\|_{2}^2}\,\mu(dx) <\infty, \ \text{for some }\varepsilon>0.
\end{equation}
The question of validating or infirming  \eqref{conjecture}  was communicated to us by Cattiaux and Guillin.

Our next goal is to disclaim \eqref{conjecture} by exhibiting a counterexample $\bar{\mu}$
on $\R$. The construction is as follows: $\bar{\mu}$ will be the image of the exponential distribution $\mu_1$ under an odd, continuous, non-decreasing and Lipschitz map $\overline{T}:\R\to\R$ which verifies $|\overline{T}(x)|\leq \sqrt{|x|}$ for all $x\in \R$ but does not satisfy the growth condition \eqref{contraction} for $\alpha(x)=x^2$, which means that
\begin{equation}\label{non-holder}
\sup_{x,y\in \R}\frac{|\overline{T}(x)-\overline{T}(y)|}{\sqrt{1+|x-y|}} =\infty.
\end{equation}
Let us take for granted the existence of such a map $\overline{T}$. The fact that it is Lipschitz then implies that $\bar{\mu}$ verifies Poincar\'e and the inequality $|\overline{T}(x)|\leq \sqrt{|x|}$ easily implies that $\int e^{\varepsilon x^2}\bar{\mu}(dx)<\infty$ for all $\varepsilon<1.$ Finally, we conclude from Theorem \ref{equivalence 2} and condition \eqref{non-holder} that $\bar{\mu}$ does not verify $\T_2$ (here we use the fact that $\overline{T}$ is actually the transport map between $\mu_1$ and $\bar{\mu}$).

Now let us construct such a map $\overline{T}$. The strategy is to wait until there is enough room under the graph of $x\mapsto \sqrt{x}$ to put a linear step with slope $1$ and range of length $n$, for each $n\in \N^*$. A possible construction is as follows: let $x_n=\frac{n(n+1)}{2}$, for all $n\in \N$ and define $\overline{T}(x)= x_{n-1}+(x-x_n^2+n)$ if $x\in [x_n^2-n, x_n^2]$ and $\overline{T}(x)=x_n$ if $x\in [x_n^2, x_{n+1}^2-(n+1)]$, for all $n\in \N^*.$ This defines $\overline{T}$ on $\R^+$ and so everywhere since $\overline{T}$ is assumed to be odd. 
This map $\overline{T}$ is clearly non-decreasing and $1$-Lipschitz and it is not difficult to check that   $|\overline{T}(x)|\leq \sqrt{|x|}$ for all $x\in \R$. Finally, $\overline{T}(x_n^2)-\overline{T}(x_n^2-n)=x_n-x_{n-1}=n$ which proves \eqref{non-holder}.

\section*{Appendix}
Usually, functional inequalities are assumed to hold ``for all functions smooth enough".
When the reference probability measure is absolutely continuous with respect to Lebesgue, this formulation makes sense. Since we allow, in this paper, probability measures to have singular parts (in particular in the examples given in Section 4), we need to clarify this condition. In our definition of Poincar\'e (and log-Sobolev), we took the class of Lipschitz functions as domain of the inequality, with
\begin{equation}\label{gradient appendix}
|\nabla f|(x)=\limsup_{y\to x}\frac{|f(y)-f(x)|}{|y-x|}
\end{equation}
in the right hand side. The following proposition establishes the equivalence between this definition and others appearing in the literature. 

\begin{prop}\label{equiv poinc}Let $\mu$ be a Borel probability measure on $\R$ with the following decomposition:
$$\mu= \mu_{ac} + \mu_{s},$$
where $\mu_{ac}$ and $\mu_{s}$ are non-negative Borel measures such that $\mu_{ac}$ is absolutely continuous with respect to Lebesgue and $\mu_{s}$ is such that there is a \emph{closed} set $C$ with $\mu_{s}(C^c)=0=\mathrm{Leb}\, (C).$\\
Let $\lambda>0$; the following are equivalent
\begin{enumerate}
\item The probability measure $\mu$ verifies
$$\lambda \mathrm{Var}_{\mu} (f)\leq \int |\nabla f|^2\,d\mu,\quad \forall f \text{ Lipschitz}.$$
\item The probability measure $\mu$ verifies
$$\lambda \mathrm{Var}_{\mu} (f)\leq \int |f'|^2\,d\mu,\quad \forall f \text{ Lipschitz and of class } \mathcal{C}^1.$$
\item The probability measure $\mu$ verifies
$$\lambda \mathrm{Var}_{\mu} (f)\leq \int |f'|^2\,d\mu_{ac},\quad \forall f \text{ Lipschitz}.$$
\end{enumerate}
The same conclusion holds for the logarithmic Sobolev inequality instead of Poincar\'e inequality.
\end{prop}
We recall that according to Rademacher theorem, Lipschitz functions are Lebesgue almost everywhere differentiable, so that the right hand side of (3) is well defined.
\proof
We do the proof in the case of Poincar\'e inequality. 
We remark that when $f$ is differentiable at $x$, then $|\nabla f|(x) = |f' (x)|$. So $(1) \Rightarrow (2)$ and $(3) \Rightarrow (1)$.\\
Let us show that (2) implies (3). First notice that (2) is equivalent to
\begin{equation}\label{eq appendix}
\lambda\mathrm{Var}_{\mu} (F_{f}) \leq \int f^2\,d\mu,
\end{equation}
for all bounded continuous $f$ and with $F_{f}(x)=\int_{0}^x f(t)\,dt.$ Take $f$ a measurable bounded function. Define $\phi_{n}(x)=\sqrt{\frac{n}{2\pi}} e^{-nx^2/2}$, $\tilde{f}_{n}=\phi_{n}\ast f$, and $h_{n}(x)=\min (1 ; n d(x,C) ),$ where $d(x,C)=\inf_{y\in C} |x-y|,$ and finally $f_{n}=\tilde{f}_{n}h_{n}.$ The functions $f_{n}$ and $\tilde{f}_{n}$ are continuous on $\R$ and it is not difficult to check that $|f_{n}|\leq |\tilde{f}_{n}| \leq M,$ where $M=\sup |f|$. Define $F_{n}=F_{f_{n}}$ and $F=F_{f}$; it holds for all $x>0$
$$|F-F_{n}|(x)\leq \int_{0}^x |\tilde{f}_{n} - f|(t)\,dt + \int_{0}^x |f|(t)(1-h_{n}(t))\,dt.$$
Since $\tilde{f}_{n} \to f$ in $\mathrm{L}_{1} ([a,b],\mathrm{Leb})$, for all bounded interval $[a,b]$,  and $1-h_{n} \to 1_{C}$ pointwise (this property requires that $C$ is closed), we easily conclude form the fact that $\mathrm{Leb}(C)=0$ that $F_{n} \to F$ pointwise. Moreover, the inequality $|F_{n}(x)|\leq M|x|$ enables to use Lebesgue dominated convergence theorem ($\mu$ has a finite moment of order $2$). So $\mathrm{Var_{\mu}}(F_{n}) \to \mathrm{Var}_{\mu}(F)$ when $n$ goes to $\infty.$ On the other hand, since $f_{n}$ is bounded and continuous, one can apply \eqref{eq appendix}, and conclude that $$\lambda\mathrm{Var}_{\mu}(F_{n})\leq \int f_{n}^2\,d\mu= \int f_{n}^2\,d\mu_{ac} \leq \int \tilde{f}_{n}^2\,d\mu_{ac},$$ where the equality follows from the fact that $f_{n}$ vanishes on $C$. It is not difficult to see that one can extract from $\tilde{f_{n}}$ a sequence converging Lebesgue almost everywhere on $\R$. Since $|\tilde{f}_{n}| \leq M$ for all $n$, one can apply Fatou's lemma along this sequence and conclude that
\begin{equation}\label{eq 2 appendix}
\lambda \mathrm{Var}_{\mu}(F)\leq \int f^2\,d\mu_{ac}, \qquad \forall f \text{ bounded.}
\end{equation}
Now, let $g$ be a Lipschitz function on $\R.$ Being Lipschitz, this function is absolutely continuous, and so its derivative $g'(t)$ exists Lebesgue almost everywhere and is in $\mathrm{L}_{1}([a,b], \mathrm{Leb})$ for all bounded interval $[a,b]$ and it holds
$$g(x)=g(0) + \int_{0}^x g'(t)\,dt,\qquad \forall x\in \R$$
(see e.g \cite{Rudin-book}).
Applying \eqref{eq 2 appendix} to the bounded function $f$ defined by $f(t)=g'(t)$ if $g$ is differentiable at $t$ and $f(t)=0$ otherwise, we finally obtain (3).
\endproof
\bibliographystyle{plain}
\bibliography{T2}

\begin{thebibliography}{10}

\bibitem{ane}
C.~An{\'e}, S.~Blach{\`e}re, D.~Chafa{\"{\i}}, P.~Foug{\`e}res, I.~Gentil,
  F.~Malrieu, C.~Roberto, and G.~Scheffer.
\newblock {\em Sur les in\'egalit\'es de {S}obolev logarithmiques}, volume~10
  of {\em Panoramas et Synth\`eses [Panoramas and Syntheses]}.
\newblock Soci\'et\'e Math\'ematique de France, Paris, 2000.
\newblock With a preface by Dominique Bakry and Michel Ledoux.

\bibitem{BCR06}
F.~Barthe, P.~Cattiaux, and C.~Roberto.
\newblock Interpolated inequalities between exponential and {G}aussian,
  {O}rlicz hypercontractivity and isoperimetry.
\newblock {\em Rev. Mat. Iberoam.}, 22(3):993--1067, 2006.

\bibitem{BK08}
F.~Barthe and A.~V. Kolesnikov.
\newblock Mass transport and variants of the logarithmic {S}obolev inequality.
\newblock {\em J. Geom. Anal.}, 18(4):921--979, 2008.

\bibitem{BR03}
F.~Barthe and C.~Roberto.
\newblock Sobolev inequalities for probability measures on the real line.
\newblock {\em Studia Math.}, 159(3):481--497, 2003.
\newblock Dedicated to Professor Aleksander Pe{\l}czy{\'n}ski on the occasion
  of his 70th birthday (Polish).

\bibitem{BGL01}
S.~G. Bobkov, I.~Gentil, and M.~Ledoux.
\newblock Hypercontractivity of {H}amilton-{J}acobi equations.
\newblock {\em J. Math. Pures Appl. (9)}, 80(7):669--696, 2001.

\bibitem{BG99bis}
S.~G. Bobkov and F.~G{\"o}tze.
\newblock Discrete isoperimetric and {P}oincar\'e-type inequalities.
\newblock {\em Probab. Theory Related Fields}, 114(2):245--277, 1999.

\bibitem{BG99}
S.~G. Bobkov and F.~G{\"o}tze.
\newblock Exponential integrability and transportation cost related to
  logarithmic {S}obolev inequalities.
\newblock {\em J. Funct. Anal.}, 163(1):1--28, 1999.

\bibitem{BH97}
S.~G. Bobkov and C.~Houdr{\'e}.
\newblock Isoperimetric constants for product probability measures.
\newblock {\em Ann. Probab.}, 25(1):184--205, 1997.

\bibitem{BH00}
S.~G. Bobkov and C.~Houdr{\'e}.
\newblock Weak dimension-free concentration of measure.
\newblock {\em Bernoulli}, 6(4):621--632, 2000.

\bibitem{BV05}
F.~Bolley and C.~Villani.
\newblock Weighted {C}sisz\'ar-{K}ullback-{P}insker inequalities and
  applications to transportation inequalities.
\newblock {\em Ann. Fac. Sci. Toulouse Math. (6)}, 14(3):331--352, 2005.

\bibitem{BS09}
A.~I. Bonciocat and K.~T. Sturm.
\newblock Mass transportation and rough curvature bounds for discrete spaces.
\newblock {\em J. Funct. Anal.}, 256(9):2944--2966, 2009.

\bibitem{CSS76}
S.~Cambanis, G.~Simons, and W.~Stout.
\newblock Inequalities for {$Ek(X,Y)$} when the marginals are fixed.
\newblock {\em Z. Wahrscheinlichkeitstheorie und Verw. Gebiete},
  36(4):285--294, 1976.

\bibitem{CG06}
P.~Cattiaux and A.~Guillin.
\newblock On quadratic transportation cost inequalities.
\newblock {\em J. Math. Pures Appl. (9)}, 86(4):341--361, 2006.

\bibitem{CGW10}
P.~Cattiaux, A.~Guillin, and L.~M. Wu.
\newblock A note on {T}alagrand's transportation inequality and logarithmic
  {S}obolev inequality.
\newblock {\em Probab. Theory Related Fields}, 148(1-2):285--304, 2010.

\bibitem{DA56}
G.~Dall'Aglio.
\newblock Sugli estremi dei momenti delle funzioni di ripartizione doppia.
\newblock {\em Ann. Scuoloa Norm. Sup. Pisa (3)}, 10:35--74, 1956.

\bibitem{Fr60}
M.~Fr{\'e}chet.
\newblock Sur les tableaux dont les marges et des bornes sont donn\'ees.
\newblock {\em Rev. Inst. Internat. Statist.}, 28:10--32, 1960.

\bibitem{Go06}
N.~Gozlan.
\newblock Integral criteria for transportation-cost inequalities.
\newblock {\em Electron. Comm. Probab.}, 11:64--77 (electronic), 2006.

\bibitem{Go07}
N.~Gozlan.
\newblock Characterization of {T}alagrand's like transportation-cost
  inequalities on the real line.
\newblock {\em J. Funct. Anal.}, 250(2):400--425, 2007.

\bibitem{Go09}
N.~Gozlan.
\newblock A characterization of dimension free concentration in terms of
  transport inequalities.
\newblock {\em Ann. Probab.}, 37(6):2480--2498, 2009.

\bibitem{Go10}
N.~Gozlan.
\newblock Poincar\'e inequalities and dimension free concentration of measure.
\newblock {\em Ann. Inst. Henri Poincar\'e Probab. Stat.}, 46(3):708--739,
  2010.

\bibitem{GL07}
N.~Gozlan and C.~L{\'e}onard.
\newblock A large deviation approach to some transportation cost inequalities.
\newblock {\em Probab. Theory Related Fields}, 139(1-2):235--283, 2007.

\bibitem{GL10}
N.~Gozlan and C.~L\'eonard.
\newblock Transport inequalities. {A} survey.
\newblock {\em Markov Process. Related Fields}, 16:635--736, 2010.

\bibitem{GRS11bis}
N.~Gozlan, C.~Roberto, and P.M. Samson.
\newblock A new characterization of {T}alagrand's transport-entropy
  inequalities and applications.
\newblock {\em Annals of Probability}, 39(3):857--880, 2011.

\bibitem{GRS12}
N.~Gozlan, C.~Roberto, and P.M. Samson.
\newblock Characterization of {T}alagrand's transport-entropy inequalities on
  metric spaces.
\newblock submitted, 2012.

\bibitem{GRS12bis}
N.~Gozlan, C.~Roberto, and P.M. Samson.
\newblock Hamilton-jacobi equations on metric spaces and transport-entropy
  inequalities.
\newblock submitted, 2012.

\bibitem{H40}
W.~Hoeffding.
\newblock Maßstabinvariante korrelationstheorie.
\newblock {\em Schriften des Mathematischen Instituts und des Instituts f\"ur
  Angewandte Mathematik der Universit\"at Berlin}, 5:181--233, 1940.

\bibitem{Ledoux-book}
M.~Ledoux.
\newblock {\em The concentration of measure phenomenon}, volume~89 of {\em
  Mathematical Surveys and Monographs}.
\newblock American Mathematical Society, Providence, RI, 2001.

\bibitem{LV09}
J.~Lott and C.~Villani.
\newblock Ricci curvature for metric-measure spaces via optimal transport.
\newblock {\em Ann. of Math. (2)}, 169(3):903--991, 2009.

\bibitem{Ma86}
K.~Marton.
\newblock A simple proof of the blowing-up lemma.
\newblock {\em IEEE Trans. Inform. Theory}, 32(3):445--446, 1986.

\bibitem{Ma91}
B.~Maurey.
\newblock Some deviation inequalities.
\newblock {\em Geom. Funct. Anal.}, 1(2):188--197, 1991.

\bibitem{Mi08}
L.~Miclo.
\newblock Quand est-ce que des bornes de {H}ardy permettent de calculer une
  constante de {P}oincar\'e exacte sur la droite?
\newblock {\em Ann. Fac. Sci. Toulouse Math. (6)}, 17(1):121--192, 2008.

\bibitem{Muc72}
B.~Muckenhoupt.
\newblock Hardy's inequality with weights.
\newblock {\em Studia Math.}, 44:31--38, 1972.
\newblock Collection of articles honoring the completion by Antoni Zygmund of
  50 years of scientific activity, I.

\bibitem{OV00}
F.~Otto and C.~Villani.
\newblock Generalization of an inequality by {T}alagrand and links with the
  logarithmic {S}obolev inequality.
\newblock {\em J. Funct. Anal.}, 173(2):361--400, 2000.

\bibitem{RR98}
S.~T. Rachev and L.~R{\"u}schendorf.
\newblock {\em Mass transportation problems. {V}ol. {I}}.
\newblock Probability and its Applications (New York). Springer-Verlag, New
  York, 1998.
\newblock Theory.

\bibitem{Ro85}
O.~S. Rothaus.
\newblock Analytic inequalities, isoperimetric inequalities and logarithmic
  {S}obolev inequalities.
\newblock 64:296--313, 1985.

\bibitem{Rudin-book}
W.~Rudin.
\newblock {\em Real and complex analysis}.
\newblock McGraw-Hill Book Co., New York, third edition, 1987.

\bibitem{St06a}
K.T. Sturm.
\newblock On the geometry of metric measure spaces. {I}.
\newblock {\em Acta Math.}, 196(1):65--131, 2006.

\bibitem{Ta96}
M.~Talagrand.
\newblock Transportation cost for {G}aussian and other product measures.
\newblock {\em Geom. Funct. Anal.}, 6(3):587--600, 1996.

\bibitem{Ta69}
G.~Talenti.
\newblock Osservazioni sopra una classe di disuguaglianze.
\newblock {\em Rend. Sem. Mat. Fis. Milano}, 39:171--185, 1969.

\bibitem{Tom69}
G.~Tomaselli.
\newblock A class of inequalities.
\newblock {\em Boll. Un. Mat. Ital. (4)}, 2:622--631, 1969.

\bibitem{Vi03}
C.~Villani.
\newblock {\em Topics in optimal transportation}, volume~58 of {\em Graduate
  Studies in Mathematics}.
\newblock American Mathematical Society, Providence, RI, 2003.

\bibitem{Vi09}
C.~Villani.
\newblock {\em Optimal transport}, volume 338 of {\em Grundlehren der
  Mathematischen Wissenschaften [Fundamental Principles of Mathematical
  Sciences]}.
\newblock Springer-Verlag, Berlin, 2009.
\newblock Old and new.

\end{thebibliography}
\end{document}